\documentclass[10pt,a4paper]{article}
\usepackage{amsmath}
\usepackage[T1]{fontenc}
\newcommand{\mathbb}{\mathbf}
\newcommand{\Pk}{{\sf P}^k(\mathbb{C})}

\newcommand{\CC}{\mathbb{C}}
\newcommand{\vol}{{\rm Vol}}
\newcommand{\htop}{h_{\rm top}}
\newcommand{\card}{{\rm Card}} 
\newcommand{\aire}{{\rm Aire}}
\newcommand{\diam}{{\rm Diam}}
\newcommand{\codim}{{\rm codim}}
\newcommand{\voir}{\emph{voir}}
\newcommand{\vect}{\wedge}
\renewcommand{\mod}{{\rm Mod}}
\newcommand{\de}{\,{\rm d}}
\renewcommand{\phi}{\varphi}
\renewcommand{\epsilon}{\varepsilon}
\newcommand{\jac}{{\rm Jac}}
\newcommand{\lov}{{\rm lov}}
\newcommand{\cC}{{\cal C}}
\renewcommand{\long}{{\rm Long}}

\DeclareFontFamily{U}{euex}{}
\DeclareFontShape{U}{euex}{m}{n}{ <-7> euex7 <8> euex8 <9> euex9 <10-> euex10}{}
\DeclareSymbolFont{eusym}{U}{euex}{m}{n}
\DeclareSymbolFontAlphabet{\eusymb}{eusym}
\DeclareMathSymbol{\intop}{\mathop}{eusym}{"52}
    \def\int{\intop\nolimits}
\DeclareMathSymbol{\produit}{\mathop}{eusym}{"51}
        \def\prod{\produit\nolimits}
\DeclareMathSymbol{\ointop}{\mathop}{eusym}{"48}
    
\DeclareMathSymbol{A}{\mathalpha}{operators}{`A}
\DeclareMathSymbol{B}{\mathalpha}{operators}{`B}
\DeclareMathSymbol{C}{\mathalpha}{operators}{`C}
\DeclareMathSymbol{D}{\mathalpha}{operators}{`D}
\DeclareMathSymbol{E}{\mathalpha}{operators}{`E}
\DeclareMathSymbol{F}{\mathalpha}{operators}{`F}
\DeclareMathSymbol{G}{\mathalpha}{operators}{`G}
\DeclareMathSymbol{H}{\mathalpha}{operators}{`H}
\DeclareMathSymbol{I}{\mathalpha}{operators}{`I}
\DeclareMathSymbol{J}{\mathalpha}{operators}{`J}
\DeclareMathSymbol{K}{\mathalpha}{operators}{`K}
\DeclareMathSymbol{L}{\mathalpha}{operators}{`L}
\DeclareMathSymbol{M}{\mathalpha}{operators}{`M}
\DeclareMathSymbol{N}{\mathalpha}{operators}{`N}
\DeclareMathSymbol{O}{\mathalpha}{operators}{`O}
\DeclareMathSymbol{P}{\mathalpha}{operators}{`P}
\DeclareMathSymbol{Q}{\mathalpha}{operators}{`Q}
\DeclareMathSymbol{R}{\mathalpha}{operators}{`R}
\DeclareMathSymbol{S}{\mathalpha}{operators}{`S}
\DeclareMathSymbol{T}{\mathalpha}{operators}{`T}
\DeclareMathSymbol{U}{\mathalpha}{operators}{`U}
\DeclareMathSymbol{V}{\mathalpha}{operators}{`V}
\DeclareMathSymbol{W}{\mathalpha}{operators}{`W}
\DeclareMathSymbol{X}{\mathalpha}{operators}{`X}
\DeclareMathSymbol{Y}{\mathalpha}{operators}{`Y}
\DeclareMathSymbol{Z}{\mathalpha}{operators}{`Z}
\DeclareMathSymbol{e}{\mathalpha}{operators}{`e}

\title{Deux caractérisations de la mesure d'équilibre d'un endomorphisme de $\Pk$}
\author{Jean-Yves Briend et Julien Duval}
\date{}
\begin{document}
\newlength{\titwidth}
\addtolength{\titwidth}{\textwidth}
\addtolength{\titwidth}{-1cm}
\newlength{\foute}
\settowidth{\foute}{\footnotesize $^a$\ }
\newlength{\addr}
\addtolength{\addr}{\titwidth}
\addtolength{\addr}{-\foute}
\noindent\begin{minipage}[t]{\titwidth}
\noindent{{\LARGE\bf Deux caractérisations de la mesure d'équilibre d'un endomorphisme de $\Pk$}}\\[.7cm]
\noindent{\bf Jean-Yves Briend$^a$ et Julien Duval$^b$}\\[.4cm]
\noindent{30 novembre 2000}\\[.7cm]
\noindent{\footnotesize $^a$\ }\begin{minipage}[t]{\addr} \footnotesize Université de Provence, Laboratoire Analyse, Topologie et Probabilités/CNRS UMR 6632, CMI, 39, rue F.~Joliot--Curie, 13453 Marseille cedex 13, France. Tél. 04 91 11 35 84. Fax 04 91 11 35 52.\\
E-mail : briend@gyptis.univ-mrs.fr\end{minipage}\\[.1cm]
\noindent{\footnotesize $^b$\ }\begin{minipage}[t]{\addr} \footnotesize Université Paul Sabatier, Laboratoire \'Emile Picard/CNRS UMR 5580, 118, route de Narbonne, 31062 Toulouse cedex 4, France. Tél. 05 61 55 63 79. Fax 05 61 55 82 00.\\
E-mail : duval@picard.ups-tlse.fr\end{minipage}\\[.4cm]
\newlength{\Resume}
\newlength{\RESUME}
\settowidth{\Resume}{{\bf\small\sf Résumé. \ }}
\addtolength{\RESUME}{\titwidth}
\addtolength{\RESUME}{-\Resume}
\newlength{\Abstract}
\newlength{\ABSTRACT}
\settowidth{\Abstract}{{\bf\small\sf Abstract. \ }}
\addtolength{\ABSTRACT}{\titwidth}
\addtolength{\ABSTRACT}{-\Abstract}

\noindent\hspace{14mm}\hrulefill \\[.1cm]
{\bf\small\sf Résumé. \ } \begin{minipage}[t]{\RESUME} \small Soit $\mu$ la mesure d'équilibre d'un endomorphisme de $\Pk$. Nous montrons ici qu'elle est son unique mesure d'entropie maximale. Nous construisons directement $\mu$ comme distribution asymptotique des préimages de tout point hors d'un ensemble exceptionnel algébrique. \end{minipage}\\[.25cm]

\nonfrenchspacing
\noindent\hspace{\Abstract} \begin{minipage}[t]{\RESUME} \bfseries\small\slshape Two characterizations of the equilibrium measure of an endomorphism of $\Pk$ \end{minipage}\\

\noindent {\bf\small\sf Abstract. \ }\begin{minipage}[t]{\RESUME} \it\small Let $\mu$ be the equilibrium measure of an endomorphism of $\Pk$. We show that it is its unique measure of maximal entropy. We build $\mu$ directly as the distribution of preimages of any point outside an algebraic exceptional set.\end
{minipage}\\[-.1cm]

\noindent\hspace{\Abstract}\hrulefill
\end{minipage}
\frenchspacing
\section*{Introduction}
John Hubbard ({\voir} [6] ainsi que [3]) a défini, pour un endomorphisme holomorphe de $\Pk$, une mesure d'entropie maximale naturelle, la \emph{mesure d'équilibre,} comme masse de Monge-Ampère d'une fonction de Green dynamique. John-Erik Forn{\ae}ss et Nessim Sibony [3] ont montré que $\mu$ était mélangeante et reflétait la distribution des préimages des points hors d'un ensemble pluripolaire. En dimension 1, depuis les travaux de Michael Ljubich [10] et Alexandre Freire, Artur Lopès et Ricardo Ma\~né [4], on peut obtenir $\mu$ directement comme distribution asymptotique des préimages hors d'un ensemble exceptionnel algébrique. Elle est, dans ce cas, l'unique mesure d'entropie maximale [10,11]. Nous montrons ici que ces méthodes s'adaptent en toute dimension avec les mêmes résultats. 

Plus précisément, soit $f : \Pk\longrightarrow\Pk$ un endomorphisme holomorphe de degré algébrique $d\geq 2$, donc de degré topologique $d^k$. On note, pour une mesure de probabilité $\mu$ sur $\Pk$, $d^{-k}f^{*}\mu$ la mesure de probabilité correspondant par dualité à la moyenne dans les fibres de $f$ pour les fonctions continues. Ainsi $\mu_{n,x}=d^{-kn}f^{n*}\delta_x$ est simplement la mesure de comptage normalisée sur $f^{-n}(x)$. L'ensemble exceptionnel $E$ de $f$ est le plus grand ensemble algébrique propre de $\Pk$ complètement invariant par $f$ (pour une application générique, $E$ est vide). C'est aussi son ensemble postcritique asymptotique : $E$ est le lieu des points $x$ tels que $\mu_{n,x}(C)$ ne tende pas vers $0$, où $C$ est le lieu critique de $f$. La première caractérisation de $\mu$ est donc le :\\[.2cm]
\noindent{\bf Théorème 1.}~---~{\it La mesure d'équilibre $\mu$ est l'unique mesure de probabilité vérifiant $d^{-k}f^*\mu=\mu$ et ne chargeant pas $E$. De plus, $\mu$ reflète la distribution des préimages des points non exceptionnels : $\mu_{n,x}\rightarrow \mu$ si et seulement si $x$ est hors de $E$.} \\[-.35cm]

La deuxième caractérisation de $\mu$ est entropique. Depuis Mikhael Gromov [5], on sait que l'entropie topologique de $f$ est $k\log d$. Par ailleurs, comme $\mu$ est de jacobien constant $d^k$, la formule de Rohlin-Parry ({\voir} [12]) nous dit que l'entropie métrique de $\mu$ vaut $k\log d$. D'après le Principe variationnel, $\mu$ est d'entropie maximale.\\[.2cm]
\noindent{\bf Théorème 2.}~---~{\it La mesure d'équilibre $\mu$ est l'unique mesure d'entropie maximale de $f$.}\\[-.3cm]

Ce résultat a été obtenu aussi par Mattias Jonsson [7] dans le cas particulier, proche de la dimension 1, des produits semi-directs de ${\sf P}^2({\bf C})$. Nos méthodes n'empruntent pas à la théorie du pluripotentiel, mais relèvent, outre de la théorie ergodique, de la géométrie algébrique et analytique élémentaire. 

Le premier résultat s'appuie sur un lemme  de construction itérative de branches inverses de $f^n$ exponentiellement contractantes, définies sur des disques plats évitant les valeurs critiques d'un itéré $f^l$ fixé, \emph{à la} Ljubich ({\voir} [10]). Quant au deuxième, il utilise la stratégie générale de Ljubich [10] en dimension 1, et une version relative de l'estimée de Gromov [5] de l'entropie topologique de $f$ par la croissance du volume de son graphe itéré.

Le texte s'organise comme suit : le lemme \emph{à la} Ljubich fait l'objet du premier paragraphe, tandis que le second s'intéresse à l'ensemble exceptionnel, le troisième achevant la preuve du théorème 1. Puis, après des rappels sur l'entropie topologique au quatrième paragraphe, le cinquième se consacre à l'unicité de la mesure d'entropie maximale. Enfin, un lemme de comparaison aire-diamètre pour les disques holomorphes, substitut au théorème de distorsion de Koebe en dimension supérieure, est détaillé dans l'appendice.

\section{Un lemme \emph{à la} Ljubich}
On construit de manière itérative, comme en dimension 1 ({\voir} [10], et aussi [4]), beaucoup de branches inverses de $f^n$ exponentiellement contractantes sur des disques plats génériques (\emph{i.e.\/} tracés sur une droite projective non contenue dans le lieu postcritique de $f$) évitant les valeurs critiques d'un itéré $f^l$ fixé. 

Notons $C$ le lieu critique de $f$, $V=f(C)$ ses valeurs critiques, $V_l=\bigcup_{q=1}^{l} f^{q}(C)$ les valeurs critiques de $f^l$; soit aussi $\tau$ le degré de $V$. Dans toute la suite, $\omega$ désignera la forme de Fubini-Study normalisée de $\Pk$. Remarquons que $f^*\omega$ est cohomologue à $d\omega$, où $d$ est le degré algébrique de $f$. \\[.2cm]
\noindent{\bf Lemme.}~---~ {\it Soit $\epsilon>0$. Il existe un entier $l\geq 0$ tel que, sur tout disque compact plat générique $\Delta$ évitant $V_l$, on peut construire $(1-\epsilon)d^{kn}$ branches inverses de $f^n$ pour $n$ assez grand, d'images $\Delta_i^{-n}$ avec $\diam(\Delta_i^{-n})\leq cd^{-n/2}$, $c$ ne dépendant pas de $n$.}\\[-.3cm]

\noindent{\bf Démonstration.} On fixe $l$ de sorte que $2\tau d^{-l}(1-d^{-1})^{-1}<\epsilon$; soit $L$ la droite projective contenant $\Delta$ ($L$ n'est, pour aucun $n$, incluse dans les valeurs critiques de $f^n$). Agrandissons légèrement $\Delta$ en $\widetilde{\Delta}$ avec les mêmes propriétés : $\widetilde{\Delta}\subset L$ et $\widetilde{\Delta}$ évite $V_l$.\\
\noindent a\/) {\it Construction des branches inverses sur $\widetilde{\Delta}$ :\/} par hypothèse, on dispose de $d^{kl}$ branches inverses de $f^l$ sur $\widetilde{\Delta}$, d'images $\widetilde{\Delta}_i^{-l}$. À l'étape suivante, une telle branche inverse de $f^l$ en engendre $d^k$ pour $f^{l+1}$ dès que $\widetilde{\Delta}_i^{-l}$ évite les valeurs critiques $V$ de $f$. Or les disques $\widetilde{\Delta}_i^{-l}$ sont disjoints, et tracés sur la courbe $f^{-l}(L)$, qui rencontre $V$ en au plus $\tau d^{(k-1)l}$ points par le théorème de Bezout. Ainsi $d^{kl}(1-\tau d^{-l})$ disques $\widetilde{\Delta}_i^{-l}$ évitent $V$ et créent $d^{k(l+1)}(1-\tau d^{-l})$ branches inverses de $f^{l+1}$ sur $\widetilde{\Delta}$. Inductivement, on obtient donc $d^{kn}(1-\tau d^{-l}(1+d^{-1}+\cdots +d^{-n+l+1}))\geq d^{kn}(1-\epsilon/2)$ branches inverses de $f^n$ sur $\widetilde{\Delta}$, d'images $\widetilde{\Delta}_i^{-n}$. \\
\noindent b\/) {\it Estimation du diamètre de leurs  images :\/} commençons par estimer l'aire de la majeure partie d'entre elles. L'aire totale de $f^{-n}(L)$ est 
\[
\aire(f^{-n}(L))=\int_{f^{-n}(L)} \omega = d^{(k-1)n},
\]
donc, comme les disques $\widetilde{\Delta}_i^{-n}$ sont disjoints sur $f^{-n}(L)$, au plus $(\epsilon/2) d^{nk}$ d'entre eux auront une aire excédant $(2/{\epsilon})d^{-n}$. Autrement dit, pour $(1-\epsilon)d^{nk}$ d'entre les disques $\widetilde{\Delta}_i^{-n}$, on aura :
\[
\aire(\widetilde{\Delta}_i^{-n})\leq \frac{2}{\epsilon}d^{-n}.
\]
Soient maintenant $\Delta_i^{-n}$ les disques plus petits obtenus en restreignant à $\Delta$ les branches inverses correspondantes : $\Delta_i^{-n}=\widetilde{\Delta}_i^{-n}\cap f^{-n}(\Delta)$. L'estimée de diamètre des $\Delta_i^{-n}$ va résulter du fait suivant (\emph{cf.\/} appendice), puisque l'anneau $\widetilde{\Delta}_i^{-n}-\Delta_i^{-n}$ a un module fixe, celui de $\widetilde{\Delta}-\Delta$ :\\
\noindent{\bf Fait.}~---~ Il existe $a>0$ tel que, pour toute paire de disques holomorphes $D\subset\widetilde{D}$ dans $\Pk$, on ait
\[
(\diam(D))^2\leq a \frac{\aire(\widetilde{D})}{\mod(\widetilde{D}-D)},
\]
ce qui achève la démonstration du lemme.\\

En conséquence, beaucoup de points de $\Pk$ ont la même distribution de préimages. Rappelons que $C$ désigne le lieu critique de $f$.\\[.2cm]
\noindent{\bf Corollaire.}~---~ {\it Soient $x,y\in \Pk$ tels que $\mu_{n,x}(C)$ et $\mu_{n,y}(C)$ tendent vers $0$ (si par exemple $x$ et $y$ sont hors du lieu postcritique de $f$). Alors $\mu_{n,x}-\mu_{n,y}$ converge faiblement vers $0$.}\\[-.3cm]

\noindent{\bf Démonstration.} L'hypothèse entraîne que $\mu_{n,x}(V_l)\rightarrow 0$ et $\mu_{n,y}(V_l)\rightarrow 0$ pour tout $l$. Fixons $\epsilon>0$, et soit $l$ donné par le lemme précédent. Alors, pour $\phi$ continue sur $\Pk$ de norme $1$, et $z,t$ hors de $V_l$, on aura
\[
\left|\int \phi\de\mu_{n,z} - \int\phi\de\mu_{n,t}\right| \leq 3\epsilon,
\]
si $n$ est assez grand. Pour voir cela, il suffit de choisir un disque $\Delta$ contenant $z$ et $t$ et évitant $V_l$ sur la droite joignant $z$ à $t$. Pour $n$ assez grand, $(1-\epsilon)d^{kn}$ branches inverses $\Delta_i^{-n}$ auront un diamètre inférieur au module de continuité de $\phi$ pour $\epsilon$, ce qui donne l'estimée. Maintenant, si $m$ est assez grand pour que $\mu_{m,x}(V_l)+\mu_{m,y}(V_l)\leq \epsilon$, on aura :
\[
\left|\int\phi\de\mu_{n,x}-\int\phi\de\mu_{n,y}\right| \leq \int\!\int\left|\int\phi\de\mu_{n-m,z}- \int\phi\de\mu_{n-m,t}\right|\de\mu_{m,x}(z)\otimes\de\mu_{m,y}(t),
\]
d'après l'identité $\mu_{n,x}=\int \mu_{n-m,z}\de\mu_{m,x}(z)$. On a donc
\[
\left|\int\phi\de\mu_{n,x}-\int\phi\de\mu_{n,y}\right| \leq 2\epsilon +\int\!\int_{V_l^c\times V_l^c}\left| \int\phi\de\mu_{n-m,z}-\int\phi\de\mu_{n-m,t}\right|\de\mu_{m,x}(z)\otimes\de\mu_{m,y}(t),
\]
et donc $\limsup_{n} |\int\phi\de\mu_{n,x}-\int\phi\de\mu_{n,y}|\leq 5\epsilon$ par le lemme de Fatou.

\section{L'ensemble exceptionnel}
C'est le lieu des points $x$ dont l'orbite négative rencontre trop souvent le lieu critique de $f$, au sens où $\mu_{n,x}(C)$ ne tend pas vers zéro. C'est donc en quelque sorte le lieu postcritique asymptotique. 

Sa description précise passe par la stratification de $\Pk$ par le degré local (ou multiplicité) de $f$. Le \emph{degré topologique local} $\deg_x f$ de $f$ en $x$ est le nombre de points de $f^{-1}(y)$ proches de $x$ pour $y$ proche de $f(x)$. Il varie entre $1$ et $d^k$ et ses strates $\{\deg_x f\geq s\}$ sont algébriques. On peut donc définir le \emph{degré topologique} $\deg_A f$ de $f$ \emph{le long d'un ensemble algébrique irréductible} $A$ par $\deg_A f=\min_{x\in A}\deg_x f$. Ce sera aussi le degré local de $f$ aux points génériques de $A$. Ce degré le long de $A$ est contr\^olé par la codimension de $A$ :\\[.2cm]
\noindent{\bf Lemme.}~---~ {\it Soit $A$ un ensemble algébrique irréductible, $s=\deg_A f$ et $p=\codim (A)$. Alors $s\leq d^p$.}\\[-.3cm]

\noindent{\bf Démonstration.} Fixons un voisinage $U$ d'un point générique $x$ de $A$, au sens où $\deg_x f=s$, $x$ un point lisse et non critique pour $f_{|A}$. Soit $\pi$ une rétraction holomorphe de $U$ sur $A\cap U$. Si $P$ est un plan projectif de codimension $p$ proche du plan tangent de $f(A)$ en $f(x)$, $\pi$ induit un revêtement ramifié de degré $s$ de $f^{-1}(P)\cap U$ sur $A\cap U$. Donc, si $Q$ est un plan projectif de dimension $p$ proche du plan tangent à $\pi^{-1}(x)$ en $x$, on aura $\card(f^{-1}(P)\cap Q)\geq s$. Cependant, comme le degré de $f^{-1}(P)$ est $d^p$, on a $\card(f^{-1}(P)\cap Q)\leq d^p$ par le théorème de Bezout, ce qui donne l'inégalité recherchée.\\  

\noindent{\bf Remarque.} Si de plus $A$ est invariant par $f$, alors $s=d^p$ équivaut à la complète invariance de $A$ : en effet, celle-ci se traduit par des égalités dans les inclusions $f^{-1}(Q)\cap A\subset f^{-1}(Q\cap A)$ pour les plans génériques $Q$ de dimension $p$, elles mêmes équivalentes aux égalités des cardinaux avec multiplicité. Or, par le théorème de Bezout, $\card(Q\cap A)=\delta$, où $\delta$ est le degré de $A$, et $\card(f^{-1}(Q)\cap A)=\delta d^{k-p}$. D'où il ressort que
\[
\sum_{z\in f^{-1}(Q\cap A)} \deg_z f = \delta d^{k} \;\mbox{et}\; \sum_{z\in f^{-1}(Q)\cap A}\deg_z f = \delta d^{k-p}s,
\]
ce qui conclut l'argument.
On voit donc que les strates de degré topologique correspondant à une puissance de $d$ vont jouer un r\^ole particulier. On pose $A_p=\{x\in\Pk, \deg_x f\geq d^p\}$.\\[.2cm]
\noindent{\bf Définition.}~---~ {\it On note $E_p$ la réunion des cycles de composantes de codimension (minimale) $p$ entièrement contenus dans $A_p$. L'ensemble exceptionnel $E$ est la réunion des $E_p$ :
\[
E=\bigcup_{p=1}^{k} E_p.
\]}\\[-.3cm]

Les propriétés suivantes se déduisent directement de ce qui précède :

\noindent 1\/) l'ensemble exceptionnel d'un itéré de $f$ coïncide avec celui de $f$. De plus $E_p$ est exactement la réunion des composantes de codimension $p$ de $A_p(f^l)$ pour $l$ assez grand.

\noindent 2\/) L'ensemble exceptionnel est le plus grand ensemble algébrique propre complètement invariant.

\noindent Le lemme suivant montre que c'est bien le lieu postcritique asymptotique recherché :\\[.2cm]
\noindent{\bf Lemme.}~---~{\it Soit $x$ hors de $E$. Alors $\mu_{n,x}(C)$ tend vers $0$.}\\[-.3cm]

Ainsi, d'après le paragraphe 1, les orbites négatives des points non exceptionnels ont même distribution :\\[.2cm]
\noindent{\bf Corollaire.}~---~{\it Soient $x$ et $y$ deux points hors de $E$. Alors $\mu_{n,x}-\mu_{n,y}$ converge faiblement vers $0$.}\\[-.3cm]

\noindent{\bf Démonstration du lemme.} Pour simplifier l'exposition, nous faisons la convention que les inégalités intervenant dans cet argument sont à constante multiplicative indépendante de $n$ près. Quitte à passer à un itéré de $f$, on peut supposer que $E_p$ coïncide avec les composantes de codimension $p$ de $A_p$.

Montrons, par récurrence descendante sur $p$, que $\mu_{n,x}(A_p)$ décroît exponentiellement vite. Par hypothèse, $f^{-n}(x)$ évite $E$, donc $\mu_{n,x}(A_k)=0$ puisque $A_k=E_k$. De plus $f^{-n}(x)$ ne rencontre $A_{p-1}$ qu'en des composantes de codimension $\geq p$. Supposons que $\mu_{n,x}(A_p)\leq \lambda_p^n$ avec $\lambda_p<1$. On va estimer $\mu_{n,x}(A_{p-1})$ en dénombrant $f^{-n}(x)\cap A_{p-1}$ et en majorant la plupart des multiplicités de $f^n$ sur $f^{-n}(x)$. Tout d'abord, $\card(f^{-n}(x)\cap A_{p-1})\leq d^{n(k-p)}$, comme on le voit par le théorème de Bezout sur $f^{-n}(P)\cap A$ pour les composantes $A$ de $A_{p-1}$ de codimension $\geq p$ et les plans génériques $P$ de dimension complémentaire passant par $x$. Ensuite, on a, pour $\rho<1$ fixé,
\[
\mu_{n,x}\left(\bigcup_{0\leq q\leq n\rho} f^{-q}(A_p)\right) \leq \sum_{n(1-\rho)\leq m\leq n} \mu_{m,x}(A_p) \leq \lambda_{p}^{(1-\rho)n}.
\]
Or, en dehors de cette réunion, on a $\deg_y (f^n)\leq ((d^p-1)^{\rho}d^{k(1-\rho)})^n$. Donc 
\[
\mu_{n,x}(A_{p-1})\leq \lambda_{p}^{(1-\rho)n} + \left(\frac{(d^p-1)^{\rho}d^{k(1-\rho)}}{d^{p}}\right)^n \leq \lambda_{p-1}^n,
\]
avec $\lambda_{p-1}<1$ si $\rho$ a été choisi de sorte que $(d^p-1)^{\rho}d^{k(1-\rho)}< d^p$. À la fin de la récurrence, on obtient $\mu_{n,x}(A_1)\leq \lambda_1^{n}$ pour un $\lambda_1<1$, et l'on conclut en passant à une décroissance exponentielle de $\mu_{n,x}(C)$ comme ci-dessus.

\section{Distribution des orbites négatives}
Rappelons que pour une mesure de probabilité $\mu$, $d^{-k}f^*\mu$ désigne la mesure de probabilité obtenue en dualisant l'opération de moyenne sur les fibres de $f$ pour une fonction continue $\phi$ sur $\Pk$ :
\[
\frac{f_*\phi}{d^k}(x)=\frac{1}{d^k}\sum_{z\in f^{-1}(x)}\phi(z).
\]
C'est un opérateur continu pour la convergence faible. Montrons la première caractérisation de la mesure d'équilibre $\mu$.\\[.2cm]
\noindent{\bf Théorème 1.}~---~{\it Il existe une unique mesure de probabilité $\mu$ sur $\Pk$ satisfaisant $d^{-k}f^*\mu = \mu$ et ne chargeant pas l'ensemble exceptionnel $E$. De plus, pour toute mesure de probabilité $\nu$ ne chargeant pas $E$, on a 
\[
\frac{f^{n*}\nu}{d^{kn}}\longrightarrow \mu.
\]
En particulier, $\mu_{n,x}=d^{-kn}f^{n*}\delta_x$ converge vers $\mu$ si et seulement si $x$ est hors de $E$.}\\[-.3cm]

\noindent{\bf Démonstration.} La convergence résulte de l'existence de $\mu$, gr\^ace au corollaire du paragraphe 2, et entraîne son unicité. En effet, si $\mu$ vérifie les hypothèses du théorème et $\nu$ ne charge pas $E$, on peut écrire
\[
\mu = \int \delta_x \de\mu(x),
\]
et donc
\[
\mu = \frac{f^{n*}\mu}{d^{kn}}=\int \mu_{n,x}\de\mu(x).
\]
De même on a
\[
\frac{f^{n*}\nu}{d^{kn}}=\int \mu_{n,y}\de\nu(y),
\]
d'où l'on déduit que
\[
\mu-\frac{f^{n*}\nu}{d^{kn}} = \int\!\int_{E^c\times E^c}(\mu_{n,x}-\mu_{n,y})\de\mu(x)\otimes\de\nu(y) \rightarrow 0,
\]
d'après le paragraphe 2. 

Il nous reste maintenant à construire $\mu$, et voici une manière de le faire. On note $\Omega=\omega^k$ la forme volume de la métrique de Fubini-Study, $\mu_n = d^{-kn}f^{n*}\Omega$ et $\nu_n = n^{-1}\sum_{m=1}^{n} \mu_m$. Par continuité de l'opérateur $d^{-k}f^*$, les valeurs d'adhérence $\nu$ de la suite $\nu_n$ dans le compact des mesures de probabilité sur $\Pk$ satisfont $d^{-k}f^*\nu = \nu$. Si l'une d'entre elles ne se concentre pas sur $E$, on pose $\mu = (\nu(E^c)^{-1}{\bf 1}_{E^c})\nu$, qui est encore un point fixe de $d^{-k}f^*$ par complète invariance de $E$.

Sinon, la suite $\nu_n$ converge vers $E$. Montrons que c'est impossible. Pour cela, désignons par $\jac(f)$ le jacobien de $f$ pour la forme $\Omega$, {\it i.e.} $\jac(f)\Omega = f^*\Omega$. On note $M=\max_{\Pk} Jac(f)$ et on choisit un voisinage $U$ de $E$ assez petit pour que $\epsilon = \max_U\jac(f)$ vérifie $\sqrt{\epsilon M}<d^k$. C'est possible car $E$ est dans le lieu critique de $f$. Comme $\nu_n(U)\rightarrow 1$, on aura, pour $n$ assez grand :
\[
\frac{3n}{4}\leq \sum_{m=1}^{n}\mu_m(U) = \int \sum_{q=0}^{n-1} {\mathbb{1}}_{U}\circ f^q \de\mu_n,
\]
car $\mu_m = f^{n-m}_*\mu_n$. Soient $r_n=\sum_{q=0}^{n-1}{\mathbb{1}}_{U}\circ f^q$ le nombre de visites de $U$ dans une $n$-orbite, et $X_n$ l'ensemble des points visitant souvent $U$ : $X_n=\{x\in\Pk, r_n(x)\geq n/2\}$. On obtient alors $\int r_n\de\mu_n \geq 3n/4$ donc $\mu_n(X_n)\geq 1/2$, ce qui, dit autrement, donne :
\[
\frac{1}{2}\leq\int_{X_n} \frac{f^{n*}\Omega}{d^{kn}} = \int_{X_n}\frac{\jac(f^n)\Omega}{d^{kn}} = \int_{X_n} \frac{ \prod_{q=0}^{n-1}(\jac(f)\circ f^q)\Omega}{d^{kn}} \leq \left(\frac{\sqrt{\epsilon M}}{d^k}\right)^n,
\]
ce qui est contradictoire pour $n$ assez grand et termine notre démonstration.\\

\noindent{\bf Remarque.} On retrouve sans mal des propriétés connues de $\mu$.

\noindent 1\/) La mesure $\mu$ est mélangeante, et donc ergodique. En effet, la convergence de $\mu_{n,x}$ vers $\mu$ se traduit dualement de la manière suivante : si $\phi$ est une fonction continue sur $\Pk$, alors $d^{-kn}f^n_* \phi$ tend vers $\int \phi\de\mu$ sur $E^c$. On aura donc, si $\psi$ est une autre fonction continue :
\[
\int\phi\;\psi\!\circ\!\! f^n\de\mu = \int \phi\;\psi\!\circ\!\! f^n \frac{f^{n*}\de\mu}{d^{kn}} = \int \phi \frac{f^{n*}(\psi\de\mu)}{d^{kn}} = \int \frac{f^n_*\phi}{d^{kn}}\psi\de\mu \rightarrow \left(\int\phi\de\mu\right)\left(\int\psi\de\mu\right),
\]
par convergence dominée. C'est le mélange (argument communiqué par Vincent Guedj).

\noindent 2\/) La contraction en $d^{-n/2}$ des branches inverses de $f^n$ dans le lemme {\it à la} Ljubich redonne la minoration des exposants de Liapounoff de $f$ relativement à $\mu$ par $(\log d)/2$ ({\voir} [1]).

\noindent 3\/) La mesure $\mu$ ne charge pas les ensembles algébriques. Sinon, $\mu$ chargerait $E$. En effet, soit $A$ un ensemble algébrique irréductible de dimension minimale chargé par $\mu$. Il existe un entier $n>0$ avec $f^{-n}(A)\supset A$, sinon les préimages $f^{-n}(A)$ seraient toutes disjointes modulo des ensembles algébriques de dimension plus petite, donc de masse nulle, et on aurait
\[
\mu\left(\bigcup_{n}f^{-n}(A)\right) = \sum_{n}\mu(f^{-n}(A)) = \sum_{n}\mu(A) =\infty.
\]
De plus, du fait de l'équation fonctionnelle $d^{-k}f^*\mu = \mu$, toutes les composantes de $f^{-n}(A)$ sont chargées par $\mu$. Donc, si $B$ était une composante de $f^{-n}(A)$ différente de $A$, on aurait
\[
\mu(f^{-n}(A))\geq\mu(A\cup B)=\mu(A)+\mu(B)>\mu(A),
\]
contredisant l'invariance de $\mu$. Ainsi $A$ est complètement invariant par $f^n$ et donc $A\subset E$. 

\section{Entropie topologique et entropie métrique}
Dans ce paragraphe nous effectuons quelques rappels. Notre référence générale est ici le livre d'Anatole Katok et Boris Hasselblatt [8]. 

\subsection{Entropie topologique}
C'est le taux de croissance du nombre de $n$-orbites discernables à $\epsilon$ près :
\[
\htop(f)=\sup_{\epsilon >0} \limsup_{n}\frac{1}{n}\log\left(\max\{ \card(F), F\;(n,\epsilon)\mbox{-séparé}\}\right),
\]
où un ensemble est $(n,\epsilon)$-séparé si deux de ses points ont des $n$-orbites $\epsilon$-séparées : si $x$ et $y$ sont deux points distincts de $F$, on a $d_n(x,y)\geq \epsilon$ pour la distance dynamique
\[
d_n(x,y)=\max_{0\leq q\leq n-1} \{d(f^q(x),f^q(y))\}.
\]
On note $B_n(x,r)$ les boules dynamiques associées.

Pour un endomorphisme holomorphe $f$ de $\Pk$, on a $\htop(f)=k\log d$. En effet, d'un c\^oté une application de classe $\cC^1$ d'une variété compacte a toujours une entropie topologique minorée par le logarithme de son degré topologique (d'après Michal Misiurewicz et Feliks Przytycki, {\voir} [8]). De l'autre, Gromov [5] obtient la majoration dans notre cadre par l'estimée $\htop(f)\leq \lov(f)$, où $\lov(f)$ est le taux de croissance du volume du graphe itéré de $f$. Si $\Gamma_n=\{(x,\ldots,f^{n-1}(x)), x\in\Pk\}$ est ce graphe, on définit
\[
\lov(f)=\limsup_{n}\frac{1}{n}\log(\vol(\Gamma_n))=\limsup_n\frac{1}{n}\log\left(\int_{\Gamma_n} \omega_n^k\right),
\]
où $\omega_n$ est la forme de K\"ahler sur $\Pk^n$ induite par la forme de Fubini-Study sur chaque facteur.

Un calcul cohomologique montre que $\lov(f)=k\log d$. La majoration $\htop(f)\leq\lov(f)$, quant à elle, repose sur le théorème de Pierre Lelong ({\voir} [9]). En effet, un ensemble $(n,\epsilon)$-séparé $F$ donne, via ses $n$-orbites, un ensemble $\epsilon$-séparé $G$ dans $\Gamma_n$ pour la distance produit, qui n'est autre que $d_n$. On a donc
\[
\vol(\Gamma_n)\geq \sum_{y\in G}\vol(B_n(y,\epsilon/2)\cap\Gamma_n),
\]
puisque les boules $B_n(y,\epsilon/2)$ sont disjointes. Le théorème de Lelong fournit une minoration indépendante de $n$ et $y$ du volume de $\Gamma_n$ dans ces boules : $\vol(B_n(y,\epsilon/2)\cap\Gamma_n)\geq c$. D'où l'on déduit que
\[
\frac{1}{n}\log(\vol(\Gamma_n))\geq \frac{\log c}{n}+\frac{1}{n}\log(\max\{\card(F), F\;(n,\epsilon)\mbox{-séparé}\}),
\]
ce qui donne la majoration souhaitée. 

L'entropie topologique se localise naturellement : \emph{l'entropie topologique de $f$ relative à $X\subset \Pk$} est
\[
\htop(f|X)=\sup_{\epsilon>0}\limsup_{n}\frac{1}{n}\log(\max\{\card(F), F\;(n,\epsilon)\mbox{-séparé},\; F\subset X\}).
\]
L'argument de Gromov s'adapte à $\htop(f|X)$ : notons $\Gamma_n|X$ la restriction de $\Gamma_n$ à $X$. Si $X$ est algébrique, alors $\htop(f|X)\leq \lov(f|X)=\limsup_{n}n^{-1}\log(\vol(\Gamma_n|X))$. Dans le cas général, on a seulement 
\[
\htop(f|X)\leq \limsup_{n}\frac{1}{n}\log(\vol((\Gamma_n|X)_{\epsilon})),
\]
pour $\epsilon>0$ fixé, où $(\Gamma_n|X)_{\epsilon}$ est le $\epsilon$-voisinage de $\Gamma_n|X$ dans $\Gamma_n$. 

\subsection{Le Principe variationnel}
Ce principe relie entropie topologique et entropie métrique. Soit $\nu$ une mesure ergodique pour $f$. Nous définissons \emph{l'entropie métrique} $h_\nu(f)$ à partir du théorème de Michael Brin et Anatole Katok ({\voir} [2]) évaluant la décroissance des masses des boules dynamiques : pour $\nu$ presque tout $x$,
\[
h_{\nu}(f)=\sup_{\epsilon>0}\liminf_{n} -\frac{1}{n}\log(\nu(B_n(x,\epsilon))).
\]
Le Principe variationnel affirme alors ({\voir} [8]) que
\[
\htop(f)=\sup\{ h_{\nu}(f), \nu\;\mbox{ergodique}\}
\]
et le supremum peut être pris sur toutes les mesures de probabilité invariantes, puisque l'entropie métrique dépend de manière affine de la mesure. En voici une version relative ({\voir} [10], lemme 7.1), conséquence du théorème de Brin et Katok :\\[.2cm]
\noindent{\bf Principe variationnel relatif :} {\it soit $X$ un borélien tel que $\nu(X)>0$. Alors $h_{\nu}(f)\leq \htop(f|X)$.}\\[-.3cm]

Le Principe variationnel pose la question, centrale en théorie ergodique, de l'existence et l'unicité d'une mesure d'entropie maximale, {\it i.e.} vérifiant $h_{\nu}(f)=\htop(f)$. Dans notre cas, la mesure d'équilibre $\mu$ en est une car elle est de jacobien constant $d^k$, {\it i.e.} pour tout borélien $B$ sur lequel $f$ est injective, on a 
\[
\mu(B)=d^{-k}\mu(f(B)).
\]
Il suffit en effet, par le théorème de Brin et Katok, d'expliciter, pour $\alpha>0$, un borélien $X_\alpha$ de mesure non nulle, avec $\mu(B_n(x,\epsilon))\leq d^{-kn(1-\alpha)}$ pour $x$ dans $X_{\alpha}$ et $n$ assez grand. Considérons pour cela un voisinage $U$ de l'ensemble des valeurs critiques $V$ assez petit pour que l'on ait $\mu(U)\leq \alpha$. Soit $X_{\alpha}$ est l'ensemble des points dont la $n$-orbite visite au plus $n\alpha$ fois $U$ pour $n$ assez grand. D'après le théorème de Birkhoff, $X_{\alpha}$ est de mesure non nulle. La masse des boules dynamiques centrées sur $X_{\alpha}$ s'estime inductivement gr\^ace à la propriété de jacobien constant de $\mu$ : soient $x\in X_{\alpha}$ et $\epsilon$ assez petit. Si $f^{q+1}(x)$ est hors de $U$, $f$ réalise une injection de $B_{n-q}(f^q(x),\epsilon)$ dans $B_{n-q-1}(f^{q+1}(x),\epsilon)$, donc $\mu(B_{n-q}(f^q(x),\epsilon))\leq d^{-k}\mu(B_{n-q-1}(f^{q+1}(x),\epsilon))$, et sinon, on a toujours $\mu(B_{n-q}(f^q(x),\epsilon))\leq \mu(B_{n-q-1}(f^{q+1}(x),\epsilon))$ par invariance de $\mu$, ce qui permet de conclure.

\section{Unicité de la mesure d'entropie maximale}
Dans ce paragraphe, nous prouvons le\\[.2cm]
\noindent{\bf Théorème 2.}~---~{\it La mesure $\mu$ est l'unique mesure d'entropie maximale de $f$.}\\[-.3cm]

\noindent{\bf Démonstration.} Il suffit de le montrer pour les mesures ergodiques. Supposons donc qu'il existe une mesure ergodique $\nu\neq\mu$ d'entropie maximale $k\log d$. Tout d'abord, remarquons que $\nu$ ne charge pas l'ensemble des valeurs critiques $V$. Sinon, d'après le paragraphe 4, on aurait $k\log d\leq \htop(f|V)\leq \lov(f|V)$. Or
\[
\vol(\Gamma_n|V) = \int_{\Gamma_n|V} \omega_n^{k-1} = \sum_{\underline{i}\in\{0,\ldots,n-1\}^{k-1}} \int_V f^{i_1 *}\omega\vect\cdots\vect f^{i_{k-1} *}\omega ,
\]
par définition de $\omega_n$, et il s'ensuit, si $\tau$ est le degré de $V$, que
\begin{align*}
\vol(\Gamma_n|V) & = \tau\sum_{\underline{i}\in\{0,\ldots,n-1\}^{k-1}} d^{i_1+\cdots +i_{k-1}}\\
                 & \leq \tau n^{k-1} d^{(k-1)n},
\end{align*}
d'où $\lov(f|V)\leq (k-1)\log d$, ce qui est contradictoire. En particulier $\nu$ ne charge pas l'ensemble exceptionnel $E$ et donc, d'après le théorème 1, n'est pas point fixe de l'opérateur $d^{-k}f^*$, {\it i.e.} n'est pas de jacobien constant $d^k$. On peut alors construire un ouvert $U$ simplement connexe (mais non nécessairement connexe) évitant $V$ avec $\nu(U)=\vol(U)=1$ sur lequel les branches inverses de $f$ ne distribuent pas équitablement $\nu$. Plus précisément, on aura $f^{-1}(U)=U_1\cup\ldots\cup U_{d^{k}}$ avec $f$ bijective entre $U_j$ et $U$, et par exemple $\nu(U_1)>d^{-k}$. Soient $\sigma$ tel que $\nu(U_1)>\sigma > d^{-k}$ et $O$ un ouvert légèrement plus petit que $U_1$ avec encore $\nu(O)>\sigma$. On se fixe également $\epsilon>0$ assez petit pour que le $\epsilon$-voisinage de $O$ soit encore inclus dans $U_1$. Suivant Ljubich, la contradiction va venir du borélien $X$ des points visitant assez souvent $O$ :
\[
X=\{x\in\Pk, r_n(x)\geq n\sigma\;\mbox{pour}\; n\geq m\},
\]
où $r_n(x)=\card(\{q, 0\leq q\leq n-1, f^{q}(x)\in O\})$. Par le théorème de Birkhoff, on peut prendre $m$ assez grand pour que $\nu(X)>0$. On a alors d'après le paragraphe 4 :
\[
k\log d \leq \htop(f|X)\leq \limsup_{n}\frac{1}{n}\log(\vol((\Gamma_n|X)_{\epsilon})).
\]
On estime le volume de $(\Gamma_n|X)_{\epsilon}$ par un codage. À un ensemble de volume nul près, $\{U_1,\ldots,U_{d^k}\}$ est une partition de $\Pk$, et elle en induit naturellement une sur $\Gamma_n$ : pour $\alpha\in\{1,\ldots ,d^k\}^{n}$, on note
\[
\Gamma_n(\alpha)=\Gamma_n \cap (U_{\alpha_0}\times\cdots\times U_{\alpha_{n-1}}).
\]
Par définition de $X$, on a, à un ensemble de volume nul près, l'inclusion
\[
(\Gamma_n|X)_{\epsilon}\subset \bigcup_{\alpha\in\Sigma_n} \Gamma_n(\alpha),
\]
où $\Sigma_n$ consiste en les $n$-uplets de $\{1,\ldots,d^k\}$ contenant beaucoup de $1$ :
\[
\Sigma_n = \{(\alpha\in\{1,\ldots,d^k\}^{n}, \card(\{q, \alpha_q = 1\})\geq n\sigma\}.
\]
Remarquons que, par un lemme de dénombrement de Ljubich (lemme 7.2 de [10]), le cardinal de $\Sigma_n$ croît moins vite que $d^{kn}$ : pour $n$ assez grand, $\card(\Sigma_n)\leq (d^{k\rho})^n$ pour un certain $\rho<1$. Ainsi, on a 
\begin{align*}
\vol((\Gamma_n|X)_{\epsilon}) & \leq \sum_{\alpha\in\Sigma_n} \int_{\Gamma_n(\alpha)} \omega_n^k \\
                              & \leq \sum_{\underline{i}\in\{1,\ldots,n-1\}^k} \sum_{\alpha\in\Sigma_n} \int_{\pi(\Gamma_n(\alpha))} f^{i_1 *}\omega\vect\cdots\vect f^{i_k *}\omega,
\end{align*}
où $\pi$ est la projection de $(\Pk)^n$ sur le premier facteur. Cela suffit pour conclure en dimension $k=1$. En effet,
\[
\int_{\pi(\Gamma_n(\alpha))} f^{i*}\omega \leq \int_{{\sf P}^1({\bf C})} \omega = 1,
\]
puisque $f^i$ est injective sur $\pi(\Gamma_n(\alpha))$, $f$ l'étant sur chaque $U_j$. Il s'ensuit que
\[
\vol((\Gamma_n|X)_{\epsilon})\leq n\card(\Sigma_n) \leq n (d^{\rho})^n,
\]
d'où $\log d\leq \htop(f|X)\leq \rho \log d$, ce qui est une contradiction. 

En dimension supérieure, on doit raffiner cet argument. Soit $\lambda$ tel que $\rho <\lambda <1$. On scinde la somme sur les $\underline{i}\in\{1,\ldots ,n-1\}^k$ en deux parties, l'une sur $\{\lambda n,\ldots ,n-1\}^k$ et l'autre sur le complémentaire. Pour $\underline{i}\in\{\lambda n,\ldots ,n-1\}^k$, on a , en posant $q=[\lambda n]$ :
\begin{align*}
\int_{\pi(\Gamma_n(\alpha))} f^{i_1 *}\omega\vect\cdots\vect f^{i_k *}\omega & = \int_{\pi(\Gamma_n(\alpha))} f^{q *}\left(f^{i_1-q *}\omega\vect\cdots\vect f^{i_k-q *}\omega\right) \\
                   & \leq\int_{\Pk}\!\!\!\!\!\!\! f^{i_1-q *}\omega\vect\cdots\vect f^{i_k-q *}\omega, \;\mbox{car $f^q$ injective sur $\pi(\Gamma_n(\alpha))$} \\
                   & \leq d^k d^{k(1-\lambda)n}.
\end{align*}
Donc la première somme est majorée par $(nd)^k \card(\Sigma_n)d^{k(1-\lambda)n}$, soit $(nd)^k(d^{k(1+\rho-\lambda)})^n$. Quant à la deuxième, on l'estime globalement, car les $\pi(\Gamma_n(\alpha))$ forment une partition de $\Pk$, par :
\[
\sum_{\underline{i}\in\{0,\ldots n-1\}^k-\{\lambda n,\ldots ,n-1\}^k} \int_{\Pk} f^{i_1 *}\omega \vect \cdots\vect f^{i_k *}\omega \leq n^k (d^{k-1+\lambda})^n,
\]
puisque $i_1+\cdots +i_k\leq (k-1+\lambda)n$ dans ce cas. Au total, on obtient
\[
k\log d \leq \htop(f|X) \leq \max(k(1+\rho-\lambda), k-1+\lambda)\log d <  k\log d,
\]
si l'on prend $\rho<\lambda<1$, ce qui est la contradiction recherchée et termine la preuve du théorème 2.

\section*{Appendice : un lemme de comparaison aire-diamètre}
Nous présentons ce lemme dans $\Pk$, bien qu'il soit valable dans un cadre plus général :\\[.2cm]
\noindent {\bf Lemme.}~---~{\it Il existe $c>0$ tel que, pour toute paire de disques holomorphes $D\subset \widetilde{D}$ dans $\Pk$, on ait
\[
(\diam(D))^2\leq c\frac{\aire(\widetilde{D})}{\mod(A)},
\]
où $A$ désigne l'anneau $\widetilde{D}-D$.}\\[-.3cm]

\noindent{\bf Démonstration.} Posons $\aire(\widetilde{D})/\mod(A)=r^2$. En jouant sur $c$, on peut supposer $r$ petit. Quitte à diminuer $\widetilde{D}$, on peut prendre $\mod(A)<1$. En exprimant ce module comme longueur extrèmale, on obtient
\[
\inf\{(\long(\gamma))^2, \gamma\;\mbox{essentielle dans}\; A\} \leq \frac{\aire(A)}{\mod(A)}.
\]
On trouve donc une courbe $\gamma$ dans $\widetilde{D}$ entourant $D$ et courte pour la métrique de Fubini-Study :
\[
(\long(\gamma))^2\leq \frac{\aire(\widetilde{D})}{\mod(A)} = r^2.
\]
Estimons maintenant le diamètre du disque $D_{\gamma}$ bordé par $\gamma$ dans $\widetilde{D}$, et a fortiori donc celui de $D$. L'ingrédient pour cela est encore le théorème de Lelong [9] minorant par $\pi$ l'aire d'une courbe holomorphe passant par $0$ dans la boule unité de $\CC^n$. En le localisant, on a le\\
\noindent{\bf Fait :} il existe $\rho>0$ tel que, pour toute boule $B(x,r)$ de $\Pk$ de rayon $r\leq \rho$, et toute courbe holomorphe $C$ dans $B(x,r)$ passant par $x$, on ait $\aire(C)\geq r^2$.

Ceci entraîne $\diam(D_{\gamma})\leq 3r$ pour $r\leq \rho$, ce qui suffit pour conclure. En effet, sinon, on trouverait $x$ dans $D_{\gamma}$ à distance au moins $r$ de $\partial D_{\gamma}=\gamma$, puisque $\long(\gamma)\leq r$. Donc $C=D_{\gamma}\cap B(x,r)$ serait une courbe holomorphe (sans bord) dans $B(x,r)$ passant par $x$. Ainsi on aurait
\[
r^2=\frac{\aire(\widetilde{D})}{\mod(A)} > \aire(\widetilde{D}) \geq \aire(C) \geq r^2,
\]
ce qui est une contradiction. 

\noindent\begin{center}{\large\sc bibliographie}\end{center}

{\small\noindent [1] J-Y.~Briend, J.~Duval, Exposants de Liapounoff et distribution des points périodiques d'un endomorphisme de ${\bf CP}^k$, {\it Acta Math.,} {\bf 182} (1999), 143-157.\\
\noindent [2] M.~Brin, A.~Katok, On local entropy, in {\it Geometric dynamics, Lect. Notes in Math.,\/} {\bf 1007} Springer Verlag (1983), 30-38.\\
\noindent [3] J.E.~Forn{\ae}ss, N.~Sibony, Complex dynamics in higher dimension, in {\it Complex potential theory,\/} P.M.~Gauthier and G.Sabidussi ed.,  Kluwer Acad. Press (1995), 131-186.\\
\noindent [4] A.~Freire, A.~Lopes, R.~Ma\~né, An invariant measure for rational maps, {\it Bol. Soc. Brasil Mat.,\/} {\bf 14} (1983), 45-62.\\
\noindent [5] M.~Gromov, On the entropy of holomorphic maps, {\it manuscrit,\/} 1977.\\
\noindent [6] J.H.~Hubbard, P.~Papadopol, Superattractive fixed points in ${\bf C}^n$, {\it Indiana Univ. Math. J.,\/} {\bf 43} (1994), 321-365.\\
\noindent [7] M.~Jonsson, Ergodic properties of fibered rational maps, {\it Ark. Mat.,\/} {\bf 38} (2000), 281-317.\\
\noindent [8] A.~Katok, B.~Hasselblatt, {\it Introduction to the modern theory of dynamical systems,\/} Cambridge Univ. Press, Encycl. of Math. and its Appl., {\bf 54,} 1995.\\
\noindent [9] P.~Lelong, Propriétés métriques des variétés analytiques complexes définies par une équation, {\it Ann. Sci. \'Ecole Norm. Sup.,} {\bf 67} (1950), 393-419.\\
\noindent [10] M.Ju.~Ljubich, Entropy properties of rational endomorphisms of the Riemann sphere, {\it Ergodic Theory Dynamical Systems,\/} {\bf 3} (1983), 351-385.\\
\noindent [11] R.~Ma\~né, On the uniqueness of the maximizing measure for rational maps, {\it Bol. Soc. Brasil Mat.,\/} {\bf 14} (1983), 27-43.\\
\noindent [12] W.~Parry, {\it Entropy and generators in ergodic theory,\/} Benjamin Press, 1969.

}

\end{document}